\documentclass[11pt]{article}

\usepackage{graphicx}

\newcommand{\bea}{\begin{eqnarray}}
\newcommand{\eea}{\end{eqnarray}}
\newcommand{\be}{\begin{eqpation}}
\newcommand{\ee}{\end{equation}}
\newcommand{\tras}{^{{\mbox{\footnotesize\sc t}}}}

\newcommand\bu {{\bf{u}}}

\newcommand\bx {{\bf{x}}}

\newcommand\bX {{\bf{x}}}

\newcommand\wg {\widehat{g}}

\newcommand\werre {\widehat{r}}

\newcommand\bbe {\mbox{$\bbeta$}}

\newcommand\bphi {\mbox{\boldmath $\phi$}}
\newcommand\bSi {\mbox{\boldmath $\Sigma$}}

\newcommand\bxi {\mbox{\boldmath $\xi$}}

\newcommand\bbeta {\mbox{\boldmath $\beta$}}
\newcommand\etab {\mbox{\boldmath $\eta$}}

\newcommand\wetab {\widehat{\etab}}

\newcommand\wbeta {\widehat{\mbox{\boldmath $\beta$}}}

\newcommand\wtbeta {\widetilde{\bbeta}}

\newcommand\wphi {\widehat{\phi}}

\newcommand{\wfi} {\widehat{\phi}}

\newcommand{\wbfi} {\widehat{\mbox{\boldmath $\phi$}}}

\newcommand\wgama {\widehat{\gamma}}
\newcommand\wbgama {\widehat{\mbox{\boldmath $\gamma$} }}

\def\MSE{\mathop{\rm MSE}}
\def\real{\hbox{$\displaystyle I\hskip -3pt R$}}

\def\realito{\tiny{\real}}

\newcommand{\convpp}{ \buildrel{a.s.}\over\longrightarrow}
\newcommand{\convprob}{ \buildrel{p}\over\longrightarrow}
\newcommand{\convdist}{ \buildrel{{\cal D}}\over\longrightarrow}

\def\dst{\displaystyle}
\def\noi{\noindent}

\def\square{\ifmmode\sqr\else{$\sqr$}\fi}
\def\sqr{\vcenter{
         \hrule height.1mm
         \hbox{\vrule width.1mm height2.2mm\kern2.18mm
\vrule width.1mm}
         \hrule height.1mm}}

\usepackage{color}

\begin{document}
\title{Robust estimators in partly linear regression models on Riemannian manifolds}

\author{Guillermo Henry and Daniela Rodriguez \\
{\small \sl Facultad de Ciencias Exactas y Naturales, Universidad de Buenos Aires and
CONICET, Argentina}}

\date{}
\maketitle

\begin{abstract}
Under a partially linear models  we study a family of robust estimates for the regression parameter and the regression function when some of the predictor variables take values on a Riemannian  manifold. We obtain  the consistency and the asymptotic normality of the proposed estimators. Also, we consider a robust cross validation procedure to select the smoothing parameter. Simulations and application to real data show  the performance of our proposal under small samples and contamination.
\end{abstract}


\section{Introduction}
Partially linear regression models (PLM) assume that the regression function can be modeled linearly on
some covariates, while it depends nonparametrically on some others. To be more precise, assume
that we have a response $y_i \in \real$   and covariates $ (\bx_i,  t_i)$ such that $\bx_i \in \real^p, t_i \in [0, 1]$ satisfying
\begin{equation}
y_i=\bx^{\tras}_i \bbeta+g(t_i)+\varepsilon_i \quad\quad 1\leq i \leq n\;,
\label{modelo}
\end{equation}
where the errors $\varepsilon_{i}$ are independent and independent of $(\bx_i^{\tras},t_i)$. Since the introductory work by \cite{eng}, the partly linear models have become an important tool in the modeling of econometric or biometric data, since they combine the flexibility of nonparametric models and the simple interpretations of the linear ones.  However, in many applications, the predictors variables take values on a Riemannian manifold more than on Euclidean space and  this structure of the variables needs to be taken into account in the estimation procedure.

In a recent work (see \cite{whr}), we studied a PLM when the explanatory variables takes values on a Riemannian manifold and we  explored  the potencial of  this model in an applications of an environment problem.  Unfortunately, as we will see in Section \ref{simu}, this approach may not work as desired because  PLM can be very sensitive to the presence of a small proportion of observations that deviate from the assumed model.   One way to avoid this problem is to derive robust estimators to fit PLM models that can resist the effect of a small number of atypical observations.  The goal of this paper is to introduce resistant estimators for the regression parameter and the regression function under PLM (\ref{modelo}), when the predictor variable $t$ takes values on a Riemannian manifolds.

This paper is organized as follows. In Section \ref{estimadores}, we give a brief summary of the  classical proposal  of estimation for this model and we introduce the robust estimates. In Section \ref{asymp}, we  study the consistency and the asymptotic distribution of the regression parameter under regular assumptions on the bandwidth sequence. A robust cross validation method for the bandwidth selection is considered in Section \ref{ventana}. Section \ref{simu} include a simulation study in order to explore the performance of the new estimators under normality and contamination. We show an example using real data, in Section \ref{real}. Proofs are given in the Appendix.

\section{The model and the estimators} \label{estimadores}
\subsection{Classical estimators}
Assume that we have a sample of $n$ independent variables   $(y_i,\bx_i^{\tras},t_i)$ in $\real^{p+1} \times M$ with identically  distribution to $(y,\bx^{\tras},t)$, where $(M,\gamma)$ is a Riemannian manifold of dimension $d$. The partially linear model assume that the relation between the response variable $y_i$ and the covariates $(\bx_i^{\tras},t_i)$ can be represented  as
\begin{equation}
y_i=\bx^{\tras}_i \bbeta+g(t_i)+\varepsilon_i \quad\quad 1\leq i \leq n\;,
\label{semipara}
\end{equation}
where the errors $\varepsilon_{i}$ are independent and independent of $(\bx_i^{\tras},t_i)$ and we will assume that $\varepsilon$ has symmetric distribution.   Denote $\phi_0(\tau)=E(y|t=\tau)$ and $\bphi(t)=(\phi_1(t),\dots,\phi_p(t))$ where $\phi_j(\tau)=E(x_{ij}|t=\tau)$ for $1\leq j\leq p$,  then we have that $g(t)=\phi_0(t)-\bphi(t)^{\tras}\bbeta$ and hence, $y-\phi_0(t)=(\bx-\bphi(t))^{\tras}\bbeta+\varepsilon$.  The classical least square estimator of $\bbe$, $\wbeta_{ls}$ can be obtained by minimizing
\begin{eqnarray*}
\widehat{\bbeta}_{ls}=\mbox{arg} \min_{\bbeta} \sum_{i=1}^n[(y_i-\wfi_{0,ls}(t_i))-(\bX_i-\wbfi_{ls}(t_i))^{\tras}\bbeta]^2,
\end{eqnarray*}
with $\wfi_{0,ls}$  and $\wbfi_{ls}$ are nonparametric kernel estimators of $\phi_0$ and $\bphi(t)$, respectively. More precisely, the nonparametric estimators  $\wfi_{0,ls}$ and $\wfi_{j,ls}$ of $\phi_0$ and $\phi_{j}$ can be defined as (see \cite{bp}),
\begin{equation}\label{estimadornop}
\wfi_{0,ls}(t)=\sum_{i=1}^n w_{n,h}(t,t_i) y_i \quad \mbox { and } \quad  \wfi_{j,ls}(t)=\sum_{i=1}^n w_{n,h}(t,t_i) x_{ij}
\end{equation}
where $w_{n,h}(t,t_i)={{\theta^{-1}_t(t_i)}K(d_{\gamma}(t,t_i)/h_n)}/[{\sum_{k=1}^n {\theta^{-1}_t(t_k)}K(d_{\gamma}(t,t_k)/h_n)}]^{-1}$
with $K:\real \to \real$  a non-negative function, $d_{\gamma}$ the distance induced by the metric ${\gamma}$,  $\theta_t(s)$  the volume density function on $(M,\gamma)$ and the bandwidth $h_n$  is a sequence of real positive numbers such that $\lim_{n\to \infty}h_n=0$ and $h_n$ are smaller than the injectivity radius of $(M,\gamma)$ ($inj_{\gamma} M$ ). As in \cite{hrnp} we consider $({M},{\gamma})$ a $d-$dimensional compact oriented Riemannian manifold without boundary. Note that in this case $inj_{\gamma} M>0$ .
 For a rigorous definition of the volume density  function and the injectivity radius  see \cite{Besse} or \cite{hriv}.

The final least square estimator of $g$ can be taken as $\wg_{ls}(t)=\wfi_{0,ls}(t)-\wbfi_{ls}(t)^{\tras}\wbeta_{ls}.$ The properties of these estimators have been studied in \cite{whr} and in the case of Euclidean data  have been widely studied in the literature, see for example \cite{eng}, \cite{sp}, \cite{chen} and \cite{liang}.

\subsection{Robust estimates}
As in the Euclidean setting, the estimators introduced by \cite{bp} are a weighted average of the response variables, these estimates are very sensitive to large fluctuations of the variables and so, the final estimator of $\bbeta$ can be seriously affected by anomalous data, as mentioned in the Introduction. To overcome this problem, \cite{hrnp} considered two families of robust estimators for the regression function when the explanatory variables $t_i$ take values on a Riemannian manifold $(M, \gamma)$. The first family combines the ideas of robust smoothing in Euclidean spaces with the
kernel weights introduced in \cite{bp}. The second generalizes to our setting the proposal given by \cite{bofr}, who considered robust nonparametric estimates using nearest neighbor weights when the predictors $t$ are on $\real^d$.

Based on the robust nonparametric estimators proposed in \cite{hrnp} and the ideas considered in \cite{bobi} for the partially linear models in the Euclidean cases, we proposed a class of estimates based on a three-step robust procedure  under the partly linear model  when some of the predictors takes values on a Riemannian manifolds. The  three-step robust estimators are defined as follows:
\begin{enumerate}
\item[] \bf Step 1: \rm  Estimate  $\phi_j(t)$, $0\leq j \leq p$ through a robust smoothing. Denote by  $\wfi_{j,{\mbox{\sc r}}}$ the obtained estimates and  $\wbfi_{\mbox{\sc r}}(t)=(\wfi_{1,{\mbox{\sc r}}}(t),\dots,\wfi_{p,{\mbox{\sc r}}}(t))^{\tras}$.
\item[]\bf Step 2: \rm  Estimate the regression parameter by applying a robust regression estimate to the residuals $y_i-\wfi_{0,{\mbox{\sc r}}}(t_i)$ and $\bx_i-\wbfi_{\mbox{\sc r}}(t_i)$. Denote by $\wbeta_{\mbox{\sc r}}$ the obtained estimator.
\item[]\bf Step 3: \rm  Define the robust estimate of the regression function $g$ as $\wg_{\mbox{\sc r}}(t)=\wfi_{0,{\mbox{\sc r}}}(t)-\wbeta_{\mbox{\sc r}}^{\tras}\wbfi_{\mbox{\sc r}}(t)$.
\end{enumerate}
Note that in the \bf Step 1\rm , the regression functions correspond to  predictors taking values in a Riemannian manifold.  Local $M-$type estimates $\wfi_{0,{\mbox{\sc r}}}$ and $\wfi_{j,{\mbox{\sc r}}}$ are defined in  \cite{hrnp} as the solution of
\begin{equation}\label{estimador}
\sum_{i=1}^n
w_{n,h}(t,t_i)\Psi\left(\frac{y_i-\wfi_{0,{\mbox{\sc r}}}(t)}{\sigma_{0,n}(t)}\right)=0 \quad \mbox{ and }\quad \sum_{i=1}^n
w_{n,h}(t,t_i)\Psi\left(\frac{x_{ij}-\wfi_{j,{\mbox{\sc r}}}(t)}{\sigma_{j,n}(t)}\right)=0
\end{equation}
respectively, where the score function $\Psi$ is strictly increasing, bounded and continuous and $\sigma_{0,n}(\tau)$ and  $\sigma_{j,n}(\tau)$ $1\leq j\leq p $  are local robust estimates.

Possible choice for the score function $\Psi$ can be the Huber or  the bisquare $\Psi$-function.
The local robust scale estimates $\sigma_{0,n}(\tau)$ and  $\sigma_{j,n}(\tau)$ $1\leq j\leq p $ can be taken as the local median of the absolute
deviations from the local median (local $\mbox{\sf MAD}$), \sl i.e. \rm the $\mbox{\sf MAD}$ (see \cite{hu}) with respect to the distributions
\begin{equation}\label{ecuacion}
F_n(y|t=\tau)=\sum_{i=1}^n w_{n,h}(\tau,t_i) I_{(-\infty,y]}(y_i) \quad\mbox{ and } \quad F_{j,n}(x|t=\tau)=\sum_{i=1}^n w_{n,h}(\tau,t_i) I_{(-\infty,x]}(x_{ij}).
\end{equation}
respectively.

In the Step 2, the robust estimation of the regression parameter can be performed by applying to the residuals any of the robust methods proposed for linear regression.
For example, we can consider M-estimates (\cite{hu}) and GM-estimators (\cite{mall}). On the other hand, high breakdown point estimates with high eficiency as MM-estimates could be evaluated (\cite{yo} and \cite{yoza}).

We consider $\wbeta_{\mbox{\sc r}}$  the solution of
\begin{equation}
\sum_{i=1}^n \psi_1\left((\werre_i-\wetab_i^{\tras}\wbeta_{\mbox{\sc r}})/{s_n}\right)  w_1\left(\|\wetab_i\|\right)\wetab_i = 0,
\label{gm}
\end{equation}
with $s_n$ a robust consistent estimate of $\sigma_{\varepsilon}$, $\werre_i=y_i-\wfi_{0,{\mbox{\sc r}}}(t_i)$, $\wetab_i=\bx_i-\wbfi_{\mbox{\sc r}}(t_i)$,
$\psi_1$ a score function and $w_1$ a weight function.  The zero of this equation can be computed iteratively using reweighting, as described for the location setting in [\cite{mmy}, Chapter 2].

The estimator defined by \cite{bp} corresponds to the choice $\Psi(u)=u$ with the estimators of the conditional distribution based on kernel weights defined in (\ref{ecuacion}). Therefore, if we considered the least square estimators of $\bbeta$ in the \bf Step 2 \rm, we obtain  the classical estimators  proposed in \cite{whr}. On the other hand,  when $(M,\gamma)$ is $ \real^d$ endowed with the canonical metric, the estimation procedure  reduces to proposal  introduced in \cite{bobi}. Details over the procedure to computing the robust nonparametric estimators in the \bf Step 1 \rm can be found in \cite{hrnp}.

\section{Asymptotic results}\label{asymp}
The theorems of this Section study the asymptotic behavior of the regression parameter estimator of the model under standard  conditions.
Let $U$ be an open set of $M$, we denote by $C^k(U)$ the set of $k$ times continuously differentiable
functions from $U$ to $\real$. As in \cite{bp}, we assume that the image measure of $P$ by $t$ is absolutely continuous with respect to the Riemannian volume measure $\nu_\gamma$ and we denote by $f$  its density  on $M$  with respect to  $\nu_\gamma$.

Let $\sigma_0(\tau)$ and $\sigma_j(\tau)$ for $1\leq j\leq n$ be the  ${\mbox{\sf mad}}$ of the conditional distribution of $y_1|t=\tau$ for $j=0$ and $x_{1j}|t=\tau$ for $1 \leq j\leq n$.
\subsection{Consistency}\label{consist}
To derive strong consistency result of the estimate $\wbeta_{\mbox{\sc r}}$ of $\bbeta$ defined in \bf Step 2 \rm, we will consider the following set of assumptions.
\begin{enumerate}
\item[$H1.\label{H1}$]
 $\Psi: \real \to \real$ is an odd, strictly increasing, bounded and continuously differentiable function, such that $ u\Psi'(u)\leq \Psi(u)$ for $u>0$.
\item[$H2.$] $F(y|t=\tau)$ and $F_j(x|t=\tau)$  are  symmetric around $\phi_0(\tau)$ and $\phi_j(\tau)$ and there are  continuous functions of $y$ and $x$ for each $\tau$.
\item[$H3.$] ${M}_0$ is a compact set on ${M}$ such that:
\begin{enumerate}
\item[i)] The density function $f$ of $t$, is a bounded function such that $\inf_{\tau\in {M}_0}f(\tau)=A>0$.

\item[ii)] $\dst\inf_{\stackrel{\tau\in {M}_0}{s\in {M}_0}} \theta_{\tau}(s)=B>0$.
\end{enumerate}
\item[$H4.$] The following equicontinuity condition holds
$$\forall \varepsilon >0,\quad \exists \delta>0: |z-z'|<\delta \Rightarrow\sup_{s\in {M}_0}|G_s(z)-G_s(z')|<\varepsilon\;$$
for the functions $G_s(z)$ equal to $F(z|t=s)$ and $F_j(z|t=s)$ for $1\leq j\leq p$.
\item[$H5.$] For any open set $U_0$ of $M$ such that $M_0\subset
U_0$,
\begin{enumerate}
\item[i)] $f$ is of class $C^2$ on $U_0$.
\item[ii)] $F(y|t=\tau)$ and $F_j(x|t=\tau)$ are  uniformly Lipschitz in
$U_0$, that is, there exists a constant $C>0$ such that
$|G_{\tau}(z)-G_{s}(z)|\leq C\,d_g(\tau,s)$ for all $\tau,s\in U_0$ and $z\in\real$,
for the functions $G_{s}(z)$ equal to $F(z|t=s)$ and $F_j(z|t=s)$ for $1\leq j\leq p$.
\end{enumerate}
\item[$H6.$] $K: \real \to \real$ is a bounded
nonnegative  Lipschitz function of order one, with compact
support $[0,1]$   satisfying $\int_{\realito^d} \bu
K(\|\bu\|)d\bu=\bf{0}$ and $0<\int_{\realito^d}
\|\bu\|^2K(\|\bu\|)d\bu<\infty$.
\item[$H7.$] The sequence $h_n$ is such that $h_n\to 0$ and  ${nh_n^d}/{\log n}\to \infty $ as $n\to \infty$.
\item[$H8.$]
The estimator $\sigma_{j,n}(\tau)$ of $\sigma_{j}(\tau)$  satisfy $\sigma_{j,n}(\tau)\convpp \sigma_j(\tau)$ as $n\to\infty$ for all $\tau\in{M}_0$ and $0\leq j\leq p.$

\end{enumerate}

\noi \textbf{Remark \ref{consist}.1.} Assumption $H1$ is a standard condition in a robustness framework.
  The fact that $\theta_s(s)=1$ for all $s\in M$  guarantees that $H3$  holds for a small compact neighborhood of $s$. $H4$ and $H5$ are needed in order to derive strong uniform consistency results.  Assumption $H6$ is a standard assumption when dealing with kernel estimators.  It is easy to see that Assumption $H8$ is satisfied, when we consider $\sigma_{j,n}(\tau)$  as the local median of the absolute
deviations from the local median.

\vskip0.1in
\noi \textbf{Theorem \ref{consist}.1.} \textsl{Under the hypothesis { $H1$} to { $H8$ }, we
have that
\begin{enumerate}
\item[a)] $ |\wbeta_{\mbox{\sc r}}-\bbeta|\convpp 0$.
\item[b)] $\sup_{\tau\in M_0}|\wg_{\mbox{\sc r}}(\tau)-g(\tau)|\convpp 0$.
\end{enumerate}}

\subsection{Asymptotic distribution}\label{distrib}
In this Section, we assume that in the \bf Step 2 \rm of the estimation procedure, the  choice for $\wbeta_{\mbox{\sc r}}$ is given as in (\ref{gm}).
More precisely, let $\psi_1$ be a score function and $w_1$ be a weight function, we will derive the asymptotic distribution of the regression parameter estimates $\wbeta_{\mbox{\sc r}}$ defined as a solution of
$$
\sum_{i=1}^n \psi_1\left((\werre_i-\wetab_i^{\tras}\wbeta_{\mbox{\sc r}})/{s_n}\right)  w_1\left(\|\wetab_i\|\right)\wetab_i = 0,$$
with $s_n$ a robust consistent estimate of $\sigma_{\varepsilon}$, $\werre_i=y_i-\wfi_{0,{\mbox{\sc r}}}(t_i)$, $\wetab_i=\bx_i-\wbfi_{\mbox{\sc r}}(t_i)$. Denote by $\etab_i=\bx_i- \bphi(t_i)$ and $r_i=y_i-\phi_0(t_i)$. Note that $r_i-\etab_i^{\tras}\bbeta=\varepsilon_i$.

To derive the asymptotic distribution of the regression parameter estimates, we will need the following set of assumptions.
\begin{enumerate}
\item[$A1.$]  $\psi_1$ is an odd, bounded and twice continuously differentiable function with bounded derivatives $\psi^\prime_1$ and $\psi^ {\prime\prime}_1$, such that the functions $u\psi^\prime_1(u)$ and $u\psi^{\prime\prime}_1(u)$ are bounded.
\item[$A2.$] $E(w_1(||\etab_1||)\; \etab_1|t_1)=0$, $E(w_1(||\etab_1||)\; ||\etab_1||^2)<\infty$ and $A=E(\psi^\prime_1(\varepsilon/\sigma_{\varepsilon})w_1(||\etab_1||)\; \etab_1{\etab_1}\tras)$ is non singular.
\item[$A3.$] The function $w_1(u)$ is bounded, Lipschitz of order 1. Moreover,   $\varphi(u)=w_1(u)u$ is also a bounded and continuously differentiable function with bounded derivative $\varphi^\prime(u)$ such that $u\varphi^\prime(u)$ is bounded.
\item[$A4.$] The functions $\phi_j(t)$ for $0\leq j \leq p$ are continuous with $\phi_j^\prime$ continuous in $M$.
\item[$A5.$]  $\wphi_j(t)$ the estimates of $\phi_j(t)$ for $0\leq j \leq p$ have first continuous derivatives in $M$ and
{\begin{eqnarray}
n^{1/4}\sup_{t\in M}|\wphi_j(t)-\phi_j(t)|\convprob 0, \mbox{ for } 0 \leq j\leq p, \label{hipoorden}\\
\sup_{t\in M}|\nabla\wphi_j(t)-\nabla\phi_j(t)|\convprob 0, \mbox{ for } 0 \leq j\leq p. \label{hipoderi}
\end{eqnarray}}
where  $\nabla \xi$ corresponds to the gradient of $\xi$ with $\xi \in {{\cal F}( M)}$ and  ${{\cal F}( M)}$  the class of functions $\{\xi\in {\cal C}^1(M): \|\xi\|_{\infty}\leq 1 \;\; \|\nabla \xi\|_{\infty}\leq 1\}$.
\item[$A6.$]
The estimator $s_n$ of $\sigma_{\varepsilon}$  satisfies $s_n\convprob \sigma_{\varepsilon}$ as $n\to\infty$.

\end{enumerate}

\vspace{0.1cm}
\noindent
\noi \textbf{Theorem \ref{distrib}.1.} Under the assumptions $A1$ to $A6$ we have that $$\sqrt{n}(\wbeta_{\mbox{\sc r}}-\bbeta)\convdist N(0,\sigma^2_{\varepsilon}A^{-1}\Sigma A^{-1}),$$ where $A$ is defined in $A2$ and $\Sigma=E(\psi^2_1(\varepsilon/\sigma_{\varepsilon}))E(w^2_1(||\etab_1||)\; \etab_1\etab_1\tras)$.

\rm
\vspace{0.2cm}

\noi \textbf{Remark \ref{distrib}.1.} To proof the previous result, we will need an inequality for the covering number of ${{\cal F}( M)}$.  The Appendix include some results related to the covering number on a Riemannian manifold.

\section{Bandwidth Selection}\label{ventana}

To select the smoothing parameter there exist  two commonly used approaches: $L^2$ cross-validation and plug-in
methods. However, these procedures may not be robust. Their sensitivity to anomalous data
was discussed by several authors, see for example \cite{lmw}, \cite{ws}, \cite{bfm},  \cite{cr} and  \cite{len}.
Under a nonparametric regression model with carriers in an Euclidean space for spline-based estimators, \cite{cr} introduced a robust cross-validation criterion to select the bandwidth parameter. Robust cross-validation selectors for kernel M-smoothers were considered in \cite{lmw}, \cite{ws} and \cite{len}, under a fully nonparametric regression model. In the Euclidean setting, for partially linear model, a robust cross-validation criterion  in the cases of   autoregression models was considered in \cite{bobi2}, while a robust plug-in procedure was studied in \cite{br}.
When the variables belong in a Riemannian manifold, a robust cross validation procedure was discussed in \cite{hrnp} under a fully nonparametric regression model, while a classical cross-validation procedure under a partly linear models was considered in \cite{whr}.

We included a robust cross-validation method for the choice of the bandwidth in the case of partially linear models that robustified the proposal given in \cite{whr}.
The robust cross-validation method constructs an asymptotically optimal data-driven bandwidth, and thus adaptive
data-driven estimators, by minimizing
$$RCV(h) =\sum_{i=1}^n\Psi^2((y_i-\wfi_{0,-i,h}(t_i))-(\bX_i-\wbfi_{-i,h}(t_i))^{\tras}\widetilde{\bbeta}),$$
where $\Psi$ is a bounded score function as the Huber's function,   $\wbfi_{-i,h}(t)=(\wfi_{1,-i,h}(t),\dots,\wfi_{p,-i,h}(t))$ and $\wfi_{0,-i,h}(t)$ denote the robust nonparametric estimators computed with bandwidth $h$ using all the data expect the $i-$th observation and $\widetilde{\bbeta}$ estimate the regression parameter by applying a robust regression estimate to the residuals $ y_i-\wfi_{0,-i,h}(t_i)$ and $\bX_i-\wbfi_{-i,h}(t_i)$.

The asymptotic properties of data–-driven estimators require further careful investigation and are beyond the scope of this paper.

\section{Simulation study}\label{simu}

In this section, we consider a simulation study designed to evaluate the performance of the robust procedure introduced in Section \ref{estimadores}.  The main objective of this study is to compare the behavior of the classical and robust estimators under normal samples and  contamination. We consider the cylinder  endowed with the metric induced by the canonical metric of $\real^3$. Because of the computational burden of the robust procedure, we performed 500 replications of independent samples of size $n=200$.  In the smoothing procedure, the kernel was taken as the quadratic kernel $K(t)=( {15}/{16}) (1-t^2)^2 I(|x|<1)$ and we choose the bandwidth using the robust cross validation procedure described in Section \ref{ventana} for the robust estimators and the classical cross validation described in \cite{whr} for the  classical estimators. The distance $d_{\gamma}$ and the volume density function for the cylinder were computed in \cite{hriv} and \cite{hrnp}. We considered the following model:

The variables $(y_i,x_i,t_i)$ for $1\leq i\leq n$ were generated as
$$
y_i=2\;x_i+ (t_{1i}+t_{2i}-t_{3i})^2+\varepsilon_i \quad \mbox{ and } \quad x_i=\sin(2t_{3i})+\eta_i
$$
where $t_i=(t_{1i},t_{2i},t_{3i})=(\cos(\theta_i),\sin(\theta_i),s_i)$ with the variables $\theta_i$ follow a uniform distribution in $(0,2\pi)$ and the variables $s_i$ are uniform in $(0,1)$, i.e. $t_i$ have support in the cylinder  with radius 1 and height between $(0,1)$.

The non contaminated cases that  denoted with $C_0$ corresponds to the errors $\varepsilon_i$ and $\eta_i$ are i.i.d. normal with mean $0$ and standard deviation $1$ and $0.05$, respectively. Besides, the so-called contaminations $C_1$ and $C_2$, which correspond to selecting a distribution in a neighborhood of the central normal distribution, are defined as $\varepsilon\sim 0.9N(0, 1) + 0.1N(0,25)$ and $\varepsilon\sim 0.9N(0, 1) + 0.1N(5,0.25)$, respectively. The contamination $C_1$ corresponds to
inflating the errors and thus, will affect the variance of the regression estimates.

\subsection{Simulation results}
Table \ref{simu}.1  shows the mean, standard deviations, mean square error for the regression estimates of $\beta$ and  the mean of the mean square error of the regression function $g$ over the 500 replications for the considered model. We denote with ${ls}$ and ${\mbox{\sc r}}$ the classical and robust estimators, respectively.  Figure \ref{simu}.1 shows the boxplot of the regression parameter.

\begin{center}
\begin{tabular}{ccccc}
\hline
& mean($\wbeta_{ls}$) & sd($\wbeta_{ls}$)& $\MSE(\wbeta_{ls})$& $\MSE(\wg_{ls})$\\
$C_0$ & 2.0732& 0.1445& 0.0262 &0.2396 \\
$C_1$ &  1.8789&  1.7592&  3.1095& 20.4485 \\
$C_2$ & 1.8722 & 1.7975 & 3.2475 &45.9719 \\\hline
& mean($\wbeta_{\mbox{\sc r}}$) & sd($\wbeta_{\mbox{\sc r}}$)& $\MSE(\wbeta_{\mbox{\sc r}})$& $\MSE(\wg_{\mbox{\sc r}})$\\
$C_0$ &2.0646& 0.1433& 0.0247& 0.2431 \\
$C_1$ & 2.0198& 0.2303& 0.0534& 0.4897 \\
$C_2$ &  2.0109 &0.2540& 0.0646& 1.3580 \\\hline
\end{tabular}
\end{center}
\vspace{-0.1cm}
\footnotesize Table \ref{simu}.1:  Performance of regression parameter and the regression functions under the different contaminations.\normalsize

The simulation results confirm the inadequate behavior of the classical estimators under the considered contaminations. The robust estimators of the regression function introduced in this work  showing only a small lack of efficiency under normality. In both cases,  the results obtained with the classical estimators are not reliable giving high mean square errors that  those corresponding to the robust procedure, under C1 and C2, respectively.
This extreme behavior of the classical estimators show its inadequacy when one suspects that the sample can contain outliers.

\newpage
\begin{center}
\hspace{-1cm}a)\hspace{8cm} b)
\vspace{-0.4cm}

\hspace{-1cm}  \includegraphics[scale=0.4]{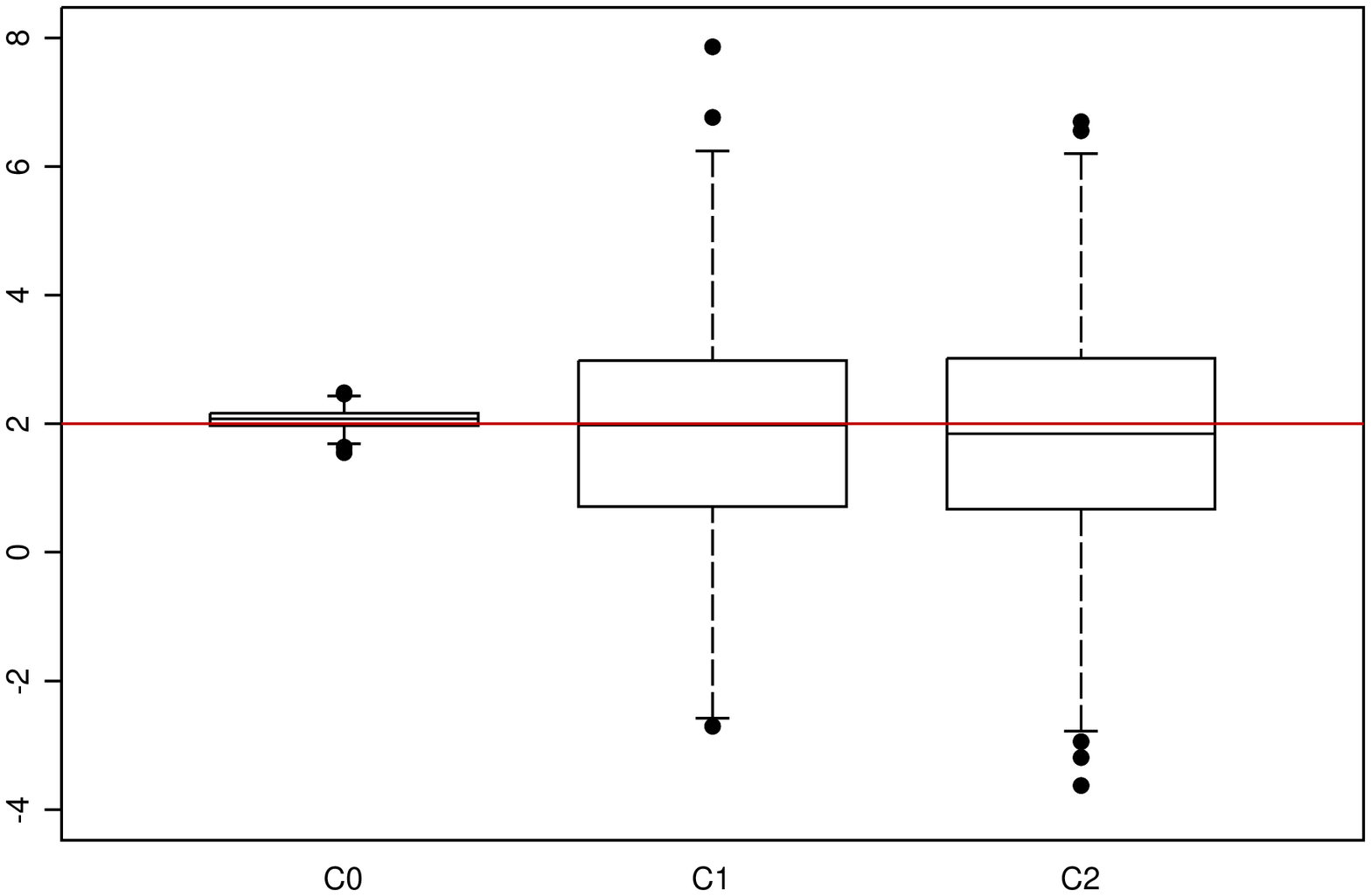}\hspace{-0.4cm}  \includegraphics[scale=0.4]{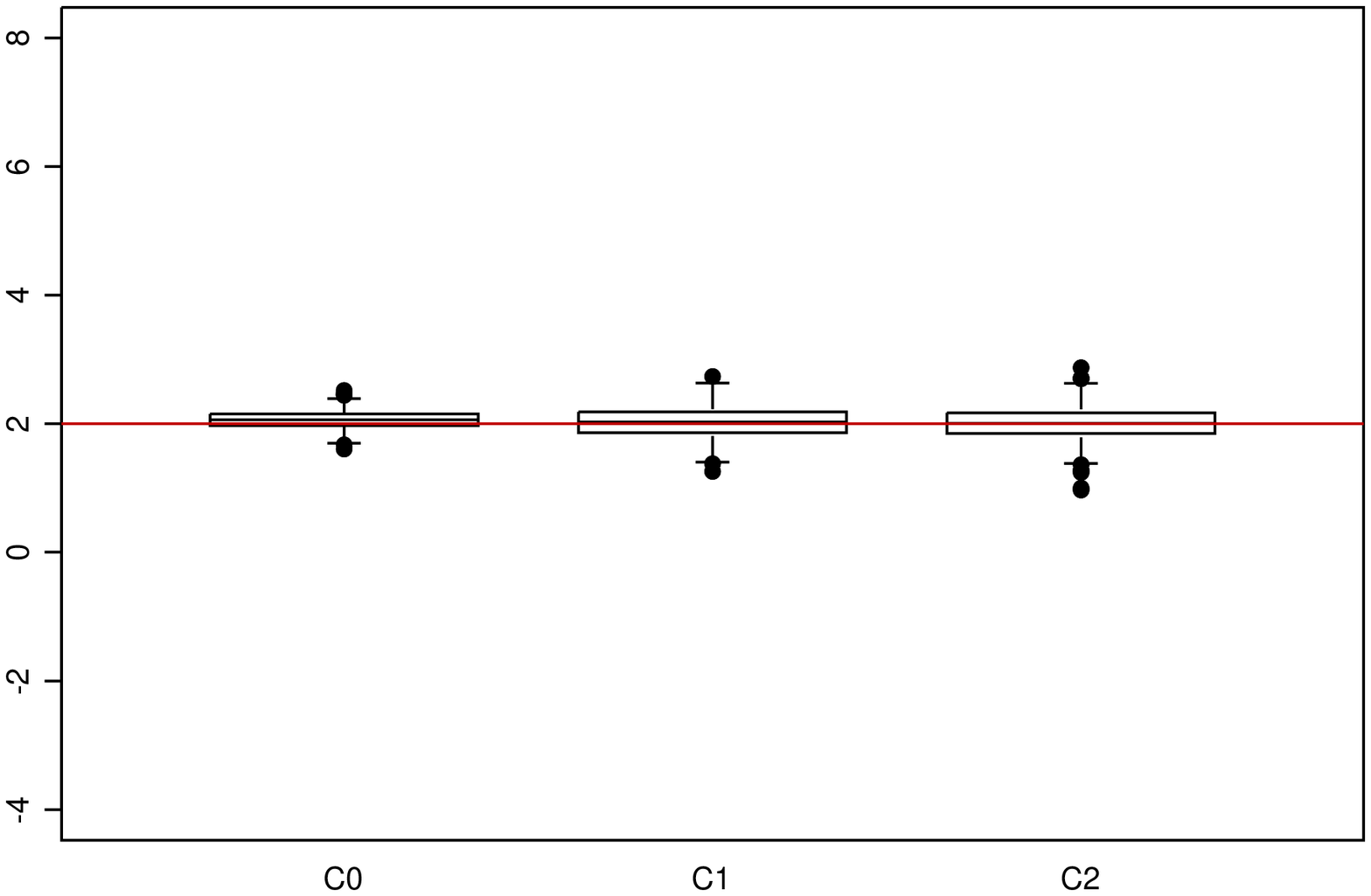}\\
\end{center}
\vspace{-0.8cm}
\footnotesize Figure \ref{simu}.1:  Boxplot of a) $\wbeta_{ls}$ the classical estimators and b) $\wbeta_{\mbox{\sc r}}$ the robust estimators under the different contaminations.\normalsize

\vspace{0.5cm}

\section{Real Example}\label{real}

The solar insolation is the amount of electromagnetic energy or solar radiation incident on the surface of the earth.  This variable measures the duration of sunlight in seconds. In the automatic stations, the World Meteorological Organization defines insolation as the sum of time intervals in which the irradiance exceeds the threshold of 120 watts per square meter. The irradiance is direct radiation normal or perpendicular to the sun on Earth's surface.  The values of the insolation in a particular location depend of the of the weather conditions and the sun's position on the horizon. For example, the presence of clouds increases the absorption, reflection and dispersion of the solar radiation. Desert areas, given the lack of clouds, have the highest values of insolation on the planet. More details  about insolation  can be seen in \cite{isola}.

As we comment above, the isolation is related with the weather conditions.  In particlar, to illustrate the proposed estimators, we will analyze the relation between the insolation, the humidity, the direction and the speed of the wind.
 We consider a data set available in http://meteo.navarra.es/. This data consists on the daily average  of relative humidity, speed and direction of the wind and the insolation. The direction's wind was measure  with the point zero in the north direction and the wind's speed was measure in meter per second. The data was measure daily in the automatic meteorologic station of Pamplona-Larrabide GN, in Navarra, Spain during the year 2004. In our study,  we consider a random sample of this dataset.

In Figure 2, we can see that the humidity and the insolation follows a lineal relation less in the outlieres contained  in the ellipse.  Therefore, we consider a partially lineal model to explain the insolation, as a linear function of the humidity and a non parametric function of the speed and direction of the wind. Note that, the variables corresponding to the wind to be modeled nonparametrically, belong to a cylinder.  In the smoothing procedure, we consider the quadratic kernel $K(t)=( {15}/{16}) (1-t^2)^2 I(|x|<1)$ and we select the bandwidth using the robust cross validation procedure for the robust estimators and the classical cross validation described in \cite{whr} for the  classical estimators.

\begin{center}
\hspace{-1cm}  \includegraphics[height=7.5cm,width=9cm]{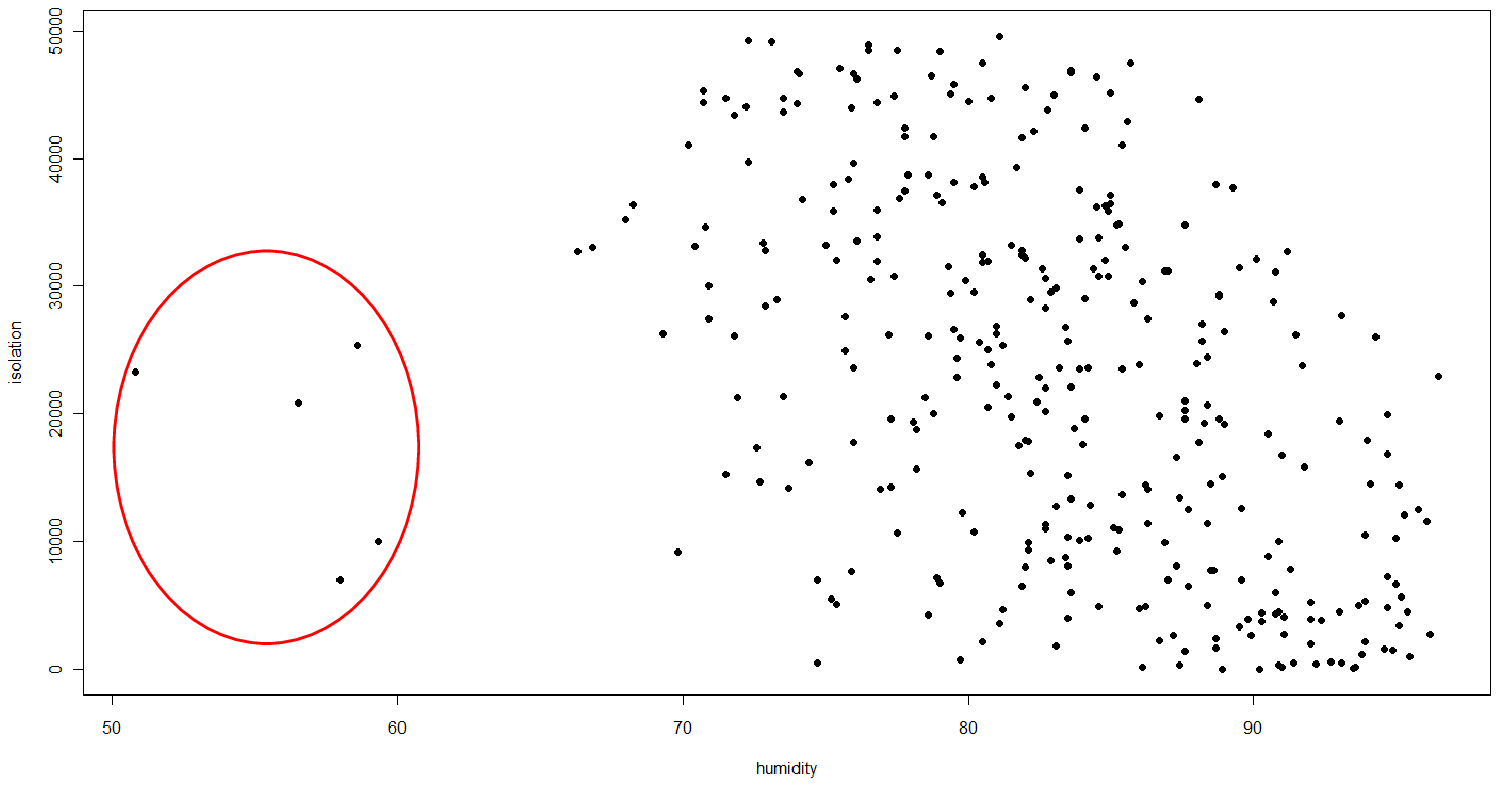}
\end{center}
\vspace{-0.5cm}
\footnotesize Figure \ref{real}.1:  Scatterplot between the insolation and humidity. The dots and the asterisks correspond to the original sample. The triangles correspond to the two outliers introduced instead of the asterisks.\normalsize
\vspace{0.5cm}

In a first step,  we apply the classical and robust methods to obtain an estimator of the regression parameter using all the data. The results were $\wbeta_{ls}=-1032.869$ and $\wbeta_{\mbox{\sc r}}=-1246.856$. Also, based in the asymptotic results obtained in Theorem 3.2.1, we calculated a confidence interval with level 0.05 in each case. To computed these confidence intervals, we   estimated  the unknown quantities. The result of the classical confidence interval was  $CI_{0.05}(\beta)= (-1229.9451  -835.7935)$ and the confidence interval based in the robust estimation was $CI_{0.05}(\beta)=( -1453.658, -1040.053)$.  On the other hand, we calculated the classical estimator using the data exept the outliers, the result was $\wbeta_{ls}=-1294.620$ and its confidence interval was $CI_{0.05}(\beta)=(-1502.983, -1086.257)$.
Thus, if we estimate the regression parameter with the classical approach when the dataset have outliers, the conclussion can be different. In the classical  case with all the data, the hypothesis that $\beta=-1000$  is rejected, while the conclusions with the classical estimator without the outliers or the robust estimators with all the data does not reject the null hypothesis.

\section*{Acknowledgments}
This research was partially supported by Grants X-018 from the Universidad de Buenos Aires, \sc pip \rm 1122008010216  from \textsc{conicet} and \sc pict \rm -00821  from \textsc{anpcyt}, Argentina.

\thispagestyle{empty} 
\appendix
\section{Appendix}\label{proofs}

\subsection{Proof of Theorem \ref{consist}.1.}

a)
Denote by $\werre_i=y_i-\wfi_{0,{\mbox{\sc r}}}(t_i)$, $\wetab_i=\bx_i-\wbfi_{\mbox{\sc r}}(t_i)$,  $\etab_i=\bx_i- \bphi(t_i)$ and $r_i=y_i-\phi_0(t_i)$. We note that $r_i=\etab_i^{\tras}\bbeta+\varepsilon_i$ and let $\widehat{P_n}(A)=\frac 1n\sum_{i=1}^n I_{A}(\werre_i,\wetab_i)$.
It is well known that the robust  regression estimates can be written as a functional of the empirical distribution. More precisely, $\wbeta =\bbeta(\widehat{P_n})$ where $\bbeta(\cdot)$ is continuous at $P$, the common distribution of $(r_i,\etab\tras_i)\tras$. Therefore, it is suffice to prove that $\Pi(\widehat{P_n},P)\convpp 0$ where $\Pi $ stands for the Prohorov distance.
Thus,  we will show that for any bounded and continuous function $f:\real^{p+1}\to\real$ we have that $|E_{\widehat{P_n}}f-E_{P}|\convpp 0$.

Note that
\begin{eqnarray*}
|E_{\widehat{P_n}}f-E_{P}f|&\leq& \frac 1n \sum_{i=1}^n|f(r_i+(\phi_0(t_i)-\wfi_0(t_i)),\etab_i+(\bphi(t_i)-\wbfi(t_i)))-f(r_i,\etab_i)|I_C(r_i,\etab_i,t_i)\\
&+&\frac 1n \sum_{i=1}^nI_{C^c}(r_i,\etab_i,t_i)
\end{eqnarray*}
where  $C_1\subset \real^{p+1}$ and $M_0\subset M$ are compact sets  such that for any $\varepsilon>0$ $P(C)> 1-\varepsilon/(4\|f\|_{\infty})$ with $C=C_1\times M_0$.

Under the assumptions by Theorem 3.3 of \cite{hrnp}, we have that $$\sup_{t\in M_0}|\wfi_{j,{\mbox{\sc r}}}(t)-\wfi_{j,{\mbox{\sc r}}}(t)|\convpp 0$$
for $0\leq j\leq p$. From this fact and the Strong Law of Large Numbers, we have that there exists a set $\aleph\subset \Omega$ such that $P(\Omega)=0$ and that for any $\omega\not\in\aleph$ we obtain that  $$\frac 1n \sum_{i=1}^nI_{C^c}(r_i,\etab_i,t_i)\to P(C^{c}).$$ Let $\bar{C_1}$ the closure of a neighborhood of radius 1 of $C_1$.The uniform continuity of $f$ on $\bar{C_1}$  implies that there exists $\delta $ such that $\max_{1\leq j\leq p+1}|u_j-u_i|$, $u,v\in \bar{ C_1}$ entails $|f(u)-f(v)|\leq \frac{\varepsilon}2$. Thus, we have that for $\omega\not\in \aleph$ and $n$ large enough $\max_{0\leq j\leq p}\sup_{t\in M_0}|\wfi_j(t)-\phi_j(t)|<\delta$ and then, for $1\leq i\leq n$, we obtain that
$$
|f(r_i+(\phi_0(t_i)-\wfi_0(t_i)),\etab_i+(\bphi(t_i)-\wbfi(t_i)))-f(r_i,\etab_i)|\leq \frac{\varepsilon}2.
$$
that conclude the proof.

b)
\square

\subsection{Entropy number}
\vskip0.1in
The main objective of this Section is to obtain an upper-bound to the entropy number of the class of functions ${{\cal F}( M)}=\{\xi\in {\cal C}^1(M): \|\xi\|_{\infty}\leq 1 \;\; \|\nabla \xi\|_{\infty}\leq 1\}$. The covering number $N(\delta,{\cal F}, \|\cdot\|)$ is the minimal number of balls, $\{\xi: \|\xi-\eta\|<\delta\}$ of radius $\delta$ needed to cover the set $\cal  F$. The entropy number is the logarithm of the covering number. This upper-bound will be use to obtain the asymptotic distribution of the regression parameter. Several authors were studied bounds to the covering  numbers for different sets, see for example \cite{van}, \cite{vv} and \cite{vw}. In particular, \cite{vw} obtained an upper-bound to the covering number for ${\cal F}(M)$ when $M$ is a bounded, convex subset of $\real^d$. For the convenience of the reader, we have included  the following remark (see \cite{ghl}).

\noi \textbf{Remark \ref{proofs}.1.}
Let  $N(\delta)$ be the minimal number of balls with radius $\delta$ needed to cover $(M,\delta)$.  A $\delta$-filling is a maximal family of pairwise disjoint open balls of radius $\delta$. We denote by $D(\delta)$, the packing number, \it i.e. the maximum number of such balls. \rm   Is easy to see that $N(2\delta)\leq D(\delta)$.
Let $diam_{(M,\gamma)}$ be the diameter of $(M,\gamma)$ and consider $\kappa\in \real$ such that $Ricc_{(M,\gamma)}\geq (d-1)\kappa$ where $Ricc_{(M\gamma)}$ is the Ricci curvature and $d$ the dimension of $M$. For example, if $\gamma$ is an Einstein metric's with scalar curvature $2(d-1)\kappa$ then the inequality is attained. Note that if $\kappa>0$ since Myers's Theorem \cite{ghl}, $(M,\gamma)$ is necessary a compact manifolds with $diam_{(M,\gamma)}\leq\pi/\sqrt{\kappa}$.   Since $M$ is compact there exists $\kappa$ with this property. Denote by $V^{\kappa}(r)$ the volume of a ball of radius $r$ in a complete, simply connected Riemannian manifold with constant curvature $\kappa$. By the Theorem of Bishop (see \cite{ghl}) we know that $\frac{Vol(B(x,r))}{V^{\kappa}(r)}$ is a non increasing  function where $B(x,r)=\{z\in M: d_{\gamma}(x,z)\leq r\}$ is the geodesic ball centered in $x$ with radius $r$. Note that, $M$ is the closure of ${B(x,diam_{(M,\gamma)})}$ for any $x\in M$. If $\{B(a_1,\frac{\delta}{2}),\dots,B(a_D,\frac{\delta}{2})\}$ with $D=D(\frac{\delta}{2})$ is a $\frac{\delta}{2}-$filling then,
$$ N(\frac{\delta}{2})\leq \frac{Vol(M)}{\inf_{1\leq i\leq D}Vol(B(a_i,\frac{\delta}{2}))}\leq \frac{V^{\kappa}(diam_{(M,\gamma)})}{V^{\kappa}(\frac{\delta}{2})}.
$$
Therefore $N(\delta)\leq C(diam_{(M,\gamma)},\kappa)\delta^{-d}$.

\noi \textbf{Lemma \ref{proofs}.1.} \textsl{Let ${{\cal F}( M)}=\{\xi\in {\cal C}^1(M): \|\xi\|_{\infty}\leq 1 \;\; \|\nabla \xi\|_{\infty}\leq 1\}$, then the covering number for the supremum norm of ${{\cal F}(M)}$ that we denote by $N(\delta,{{\cal F}(M)},\|\cdot\|_{\infty})$ satisfies that $\log N(\delta,{{\cal F}(M)},\|\cdot\|_{\infty})<A\delta^{-d}$.}

\noi\sf Proof of Lemma \ref{proofs}.1. \rm   Let ${\cal A}=\{B(a_1,\delta),\dots,B(a_N,\delta)\}$ be a covering  of $M$ by open balls of radius $\delta$. By the remark above, we may assume that $N \leq C(diam_{(M,\gamma)},\kappa)  \delta^{-d}$.  Also, we can choose the covering ${\cal A}$ such that $B(a_i,\delta)\cap B(a_{i+1},\delta)\neq \emptyset$ for $1\leq i \leq N-1 $ and $a_i\neq a_j$  for $1\leq i,j\leq N$. Let $\xi\in {{\cal F}(M)}$, we define the function $\widetilde{\xi}=\sum_{i=1}^N\delta\left[\frac{\xi(a_i)}{\delta}\right]I_{D_i}$ where $D_1=B(a_1,\delta)$, $D_i=B(a_i,\delta)\backslash \cup_{j=1}^{i-1}B(a_{j},\delta)$ and $[a]$ denotes the integer part of $a$.

Let $x\in M$ and $1\leq k\leq N$ such that $x\in D_k$, then we have that $|\widetilde{\xi}(x)-\xi(x)|\leq|\widetilde{\xi}(x)-\xi(a_k)|+|\xi(a_k)-\xi(x)|$. Since $\widetilde{\xi}(a_k)=\widetilde{\xi}(x)$ and $\xi(a_k)=\widetilde{\xi}(a_k)+\delta(\frac{\xi(a_k)}{\delta}-[\frac{\xi(a_k)}{\delta}])=\widetilde{\xi}(a_k)+\delta B$ with $0\leq B<1$, we have that $|\widetilde{\xi}(x)-\xi(x)|\leq 2\delta$.

For the first value $\widetilde{\xi}(a_1)$ of a generic function $\widetilde{\xi}$, we have $[\frac 1{\delta}]+1$ possibilities. Since,
$$
|\widetilde{\xi}(a_k)-\widetilde{\xi}(a_{k-1})|\leq |\widetilde{\xi}(a_k)-\xi(a_{k})|+|{\xi}(a_k)-{\xi}(a_{k-1})|+|{\xi}(a_{k-1})-\widetilde{\xi}(a_{k-1})|\leq 4\delta.
$$
Therefore, for each value of $\widetilde{\xi}(a_{k-1})$ we can choose $9$ possibilities for $\widetilde{\xi}(a_k)$. Then is easy to verify that
$$
N(2\delta,{{\cal F}(M)},\|\cdot\|_{\infty})\leq ([\frac 1{\delta}]+1) 9^N.
$$
which finish the proof. \square

{\noi \textbf{Remark \ref{proofs}.2.} Since $N(\delta,{{\cal F}(M)}, L^2(Q))\leq N(\delta,{{\cal F}(M)},\|\cdot\|_{\infty})$
then Lemma \ref{proofs}.1 entails that the covering number of ${{\cal F}(M)}$ satisfies, $\log N(\delta,{{\cal F}(M)}, L^2(Q))<A\delta^{-d}$.}

\subsection{ Proof of Theorem \ref{distrib}.1.}

Using a Taylor expansion around $\wbeta_{\mbox{\sc r}}$ we have that $S_n=A_n(\wbeta_{\mbox{\sc r}}-\bbeta)$ where
\begin{eqnarray*}
S_n&=&\frac{1}{n} \sum_{i=1}^n \psi_1\left((\werre_i-\wetab_i^{\tras}\bbeta)/{s_n}\right)  w_1\left(\|\wetab_i\|\right)\wetab_i \\
A_n&=&\frac{1}{n} \sum_{i=1}^n \psi^\prime_1\left((\werre_i-\wetab_i^{\tras}\wtbeta)/{s_n}\right)  w_1\left(\|\wetab_i\|\right)\wetab_i \wetab_i\tras.
\end{eqnarray*}
where $\wtbeta$ is an intermediate point between $\bbe$ and $\wbeta_{\mbox{\sc r}}$.  Analogous arguments to those used in Lemma 2 in \cite{bobi} allow to show that $A_n\convprob A$ where $A$ is defined  in $A2$.

Since $\frac{\sqrt{n}}{n} \sum_{i=1}^n \psi_1\left(\varepsilon_i/{\sigma_{\varepsilon}}\right)  w_1\left(\|\etab_i\|\right)\etab_i$ is asymptotically normally distributed with covariance $\bSi$, it will enough to show that
\begin{eqnarray}
&\sqrt{n}&[S_n-\frac{1}{n} \sum_{i=1}^n \psi_1\left(\varepsilon_i/{s_n}\right)  w_1\left(\|\etab_i\|\right)\etab_i]\convprob 0,\label{prob1}\\
&\sqrt{n}&[\frac{1}{n} \sum_{i=1}^n \psi_1\left(\varepsilon_i/{s_n}\right)  w_1\left(\|\etab_i\|\right)\etab_i-\frac{1}{n} \sum_{i=1}^n \psi_1\left(\varepsilon_i/{\sigma_{\varepsilon}}\right)  w_1\left(\|\etab_i\|\right)\etab_i]\convprob 0.\label{prob2}
\end{eqnarray}
We first prove (\ref{prob1}). Using a Taylor expansion of order two, we have that the following decomposition.
\begin{eqnarray*}
\sqrt{n}[S_n-\frac{1}{n} \sum_{i=1}^n \psi_1\left(\varepsilon_i/{s_n}\right)  w_1\left(\|\etab_i\|\right)\etab_i]&=&\sum_{i=1}^5S_{ni}
\end{eqnarray*}
where
\begin{eqnarray*}
S_{n1}&=&\frac{\sqrt{n}}{n} \sum_{i=1}^n \psi^\prime_1\left(\varepsilon_i/{s_n}\right) [\wbgama\tras(t_i)\bbeta-\wgama_0(t_i)]  w_1\left(\|\etab_i\|\right)\etab_i\\
S_{n2}&=&\frac{s_n\sqrt{n}}{n} \sum_{i=1}^n \psi_1\left(\varepsilon_i/{s_n}\right)  [w_1\left(\|\wetab_i\|\right)\wetab_i-w_1\left(\|\etab_i\|\right)\etab_i]\\
S_{n3}&=&\frac{s_n\sqrt{n}}{n} \sum_{i=1}^n [\psi_1\left(\werre_i-\wetab_i^{\tras}\bbeta/{s_n}\right)-\psi_1\left(\varepsilon_i/{s_n}\right)]  w_1\left(\|\wetab_i\|\right)[\wetab_i-\etab_i]\\
S_{n4}&=&\frac{\sqrt{n}}{2n} \sum_{i=1}^n \psi^{\prime\prime}_1\left(\varsigma_i/{s_n}\right) [\wbgama\tras(t_i)\bbeta-\wgama_0(t_i)]^2 w_1\left(\|\wetab_i\|\right)\etab_i\\
S_{n5}&=&\frac{\sqrt{n}}{n} \sum_{i=1}^n \psi_1\left(\varepsilon_i/{s_n}\right)[\wbgama\tras(t_i)\bbeta-\wgama_0(t_i)]  [w_1\left(\|\wetab_i\|\right)-w_1\left(\|\etab_i\|\right)]\etab_i\\
\end{eqnarray*}
where  $\wgama_j(t)=\wfi_j(t)-\phi_j(t)$ for $0\leq j\leq n$ and $\wbgama(t)=(\wgama_1,\dots,\wgama_n)$.
By $A3$, $A5$ and $A6$ is easy to see that $\|S_{in}\|\convprob 0$ for  $i=3,4,5.$

Let
\begin{eqnarray*}
{\cal J}^{(j)}_{1n}(\sigma,\xi)\!\!\!\!&=&\!\!\!\!\frac{\sqrt{n}}{n} \sum_{i=1}^n f^{(j)}_{1,\sigma,\xi}(r_i,\etab_i,t_i)\\
&=&\frac{\sqrt{n}}{n} \sum_{i=1}^n \psi^\prime_1\left(\frac{r_i-\etab_i\tras\bbe}{\sigma}\right) \xi(t_i) w_1\left(\|\etab_i\|\right)(\etab_i)_j\\
{\cal J}^{(j)}_{2n}(\sigma,\bxi)\!\!\!\!&=&\!\!\!\!\frac{\sqrt{n}}{n} \sum_{i=1}^n f^{(j)}_{2,\sigma,\bxi}(r_i,\etab_i,t_i)\!\! \\
&=&\!\!\frac{\sigma\sqrt{n}}{n} \sum_{i=1}^n  \psi_1\left(\frac{r_i-\etab_i\tras\bbe}{\sigma}\right)\!\!  [w_1\left(\|\etab_i+\bxi\|\right)(\etab_i+\bxi(t_i))_j-w_1\left(\|\etab_i\|\right)(\etab_i)_j]
\end{eqnarray*}
Therefore, it remains to show that ${\cal J}^{(j)}_{1n}(s_n,\wgama_s)\convprob 0$ and ${\cal J}^{(j)}_{1n}(s_n,\wbgama)\convprob 0$ for $0\leq j,s\leq p.$  From now on, we will omitted the superscript $j$ for the sake of simplicity.

Let ${{\cal F}(M)}=\{\xi\in {\cal C}^1(M): \|\xi\|_{\infty}\leq 1 \; \|\xi^\prime\|_{\infty}\leq 1\}$ and consider the classes of functions
\begin{eqnarray*}
{\cal F}_1&=&\{f_{1,\sigma,\xi}(r,\etab,t)\quad \sigma\in(\sigma_{\varepsilon}/2,2\sigma_{\varepsilon})\quad \xi\in{{\cal F}(M)}\}\\
{\cal F}_2&=&\{f_{2,\sigma,{\bxi}}(r,\etab,t)\quad \sigma\in(\sigma_{\varepsilon}/2,2\sigma_{\varepsilon})\quad \bxi=(\xi_1,\dots,\xi_p), \; \xi_s\in{{\cal F}(M)}\}
\end{eqnarray*}
Note that, the independence of $\varepsilon_i$ and $(\bx_i,t_i)$, $A2$ and the fact that the errors $\varepsilon$ have symmetric distribution imply  that $E(f(r_i,\etab_i,t_i))=0$ for any  $f\in {\cal F}_1\cup {\cal F}_2$. As in \cite{bobi}, it is easy to see that the covering number of the classes ${\cal F}_1$ and ${\cal F}_2$ satisfy
$$
N(C_1\epsilon,{\cal F}_1, L^2(Q))\leq N(\epsilon,{{\cal F}(M)},L^2(Q)) \; N(\varepsilon,(\sigma_{\varepsilon}/2,2\sigma_{\varepsilon}), |\cdot |)
$$
$$ N(C_2\epsilon,{\cal F}_2, L^2(Q))\leq N^p(\epsilon,{{\cal F}(M)},L^2(Q)) \; N(\varepsilon,(\sigma_{\varepsilon}/2,2\sigma_{\varepsilon}), |\cdot |)$$
where $Q$ is any probability measure. Since  Remark \ref{proofs}.2, the covering number of ${{\cal F}(M)}$ satisfies that $\log N(\epsilon,{{\cal F}(M)}, L^2(Q))<A\varepsilon^{-d}$. Therefore, we get that these classes have finite uniform-entropy. For $0<\delta<1, $ consider the subclasses ${\cal F}_{1,\delta}$ and ${\cal F}_{2,\delta}$ of ${\cal F}_1$ and ${\cal F}_2$ respectively, defined by,
\begin{eqnarray*}
{\cal F}_{1,\delta}&=&\{f\in{\cal F}_1\quad  \xi\in{{\cal F}(M)},\; \|\xi\|_{\infty}<\delta\}\\
{\cal F}_{2,\delta}&=&\{f\in{\cal F}_2\quad  \bxi=(\xi_1,\dots,\xi_p), \; \xi_s\in{{\cal F}(M)},\;\|\xi_s\|_{\infty}<\delta\}
\end{eqnarray*}
For any $\epsilon>0$, let $0<\delta<1$ since $A5$ and $A6$ we obtain that for $n$ large enough $P(s_n\in(\sigma_{\varepsilon}/2,2\sigma_{\varepsilon}))>1-\delta/2$ and $P(\wgama_s\in{{\cal F}(M)} \mbox{ and } \|\wgama_s\|_{\infty}<\delta)>1-\delta/2$ for $0\leq s \leq p$.

Then, the maximal inequality for covering numbers entails that for $0\leq s\leq p$
\begin{eqnarray*}
P(|{\cal J}_{1n}(s_n,\wgama_s)|>\epsilon)&\leq&P(|{\cal J}_{1n}(s_n,\wgama_s)|>\epsilon; \; s_n\in(\sigma_{\varepsilon}/2,2\sigma_{\varepsilon});\;\wgama_s\in{{\cal F}(M)} \mbox{ and } \|\wgama_s\|_{\infty}<\delta )+\delta\\
&\leq& P\left(\sup_{f\in {\cal F}_{1,\delta}}\left|\frac{\sqrt{n}}{n} \sum_{i=1}^n f(r_i,\etab_i,t_i)\right |>\epsilon\right)+\delta\\
&\leq& \frac 1{\epsilon}E\left(\sup_{f\in {\cal F}_{1,\delta}}\left|\frac{\sqrt{n}}{n} \sum_{i=1}^n f(r_i,\etab_i,t_i)\right |\right)+\delta\\
&\leq&\frac 1{\epsilon}{\cal G}(\delta,{\cal F}_1)+\delta
\end{eqnarray*}
where ${\cal G}(\delta,{\cal F})=\sup_{Q}\int_0^{\delta} \sqrt{1+\log N(\varepsilon \|F\|_{Q,2},{\cal F}, L^2(Q))}d\epsilon$, then the fact that  ${\cal F}_1$ satisfies the the uniform--entropy conditions we get that $\lim_{\delta\to 0}{\cal G}(\delta,{\cal F}_1)= 0$, therefore $S_{1n}\convprob 0$. Similarly for ${\cal J}_{2n}(s_n,\wbgama)$ and the class ${\cal F}_2$ and we get that $S_{2n}\convprob 0$.

The proof of (\ref{prob2}), follows using analogous arguments that those considered in (\ref{prob1}).\square

\vspace{1cm}\noindent {\sf Departamento de Matem\'atica,  FCEyN, Universidad de Buenos Aires\\
 Ciudad Universitaria, Pabell\'on I,  Buenos Aires, C1428EHA, Argentina  \\
e-mail address, G. Henry: ghenry@dm.uba.ar  \\
e-mail address, D. Rodriguez: drodrig@dm.uba.ar}

\end{document}